\documentclass[12pt,leqno,dvips]{article}
\usepackage{graphicx}
\begin{document}
\newtheorem{theorem}{Theorem}
\newtheorem{lemma}{Lemma}
\title{A nice example of Lebesgue integration}
\author{Joseph L. Gerver \\ Department of  Mathematics \\ Rutgers University \\ Camden, NJ 08102 \\ U.S.A. \\ gerver@camden.rutgers.edu}

\maketitle

\begin{abstract}
We explore the properties of an interesting new example of a function which is Lebesgue integrable but not Riemann integrable.
\end{abstract}

\section{Introduction}

Some years ago, while I was teaching Lebesgue's theory of integration to my real analysis class, one of the students, Michael Machuzak, asked for an honest example of a function that was Lebesgue integrable but not Riemann integrable.  He pointed out that all of my examples were the characteristic functions of Cantor sets, which he said was like developing Riemann's theory of integration, and then using it only to find the areas of rectangles.

No such example came immediately to mind, and I told Machuzak that I would get back to him.  Nor could I find any examples on the shelf of analysis textbooks in my office.  To be sure, the historical archetype of a function which is Lebesgue integrable but not Riemann integrable is the derivative of Volterra's function \cite {bressoud} (pp. 89-94).  But I would have had to spend some time constructing that function in class, and I felt that a one-line question ought to have a one-line answer.  So the following week, I gave the class the function
\begin{equation} 
f(x)=\prod_{n=0}^{\infty} [\sin (2^n x)]^{2/(2n+1)^2}.
\end{equation}

Over the next few years, I came to realize that this function has a number of interesting properties, and I thought it ought to be more well known, which is my reason for writing this paper.

Figure 1 shows the graph of $f(x)$, as plotted by Maple.  However, as we shall see, there is no truly satisfactory way to picture this graph, although fig. 1 may be as good as any.

\begin{figure}[htbp]
\includegraphics{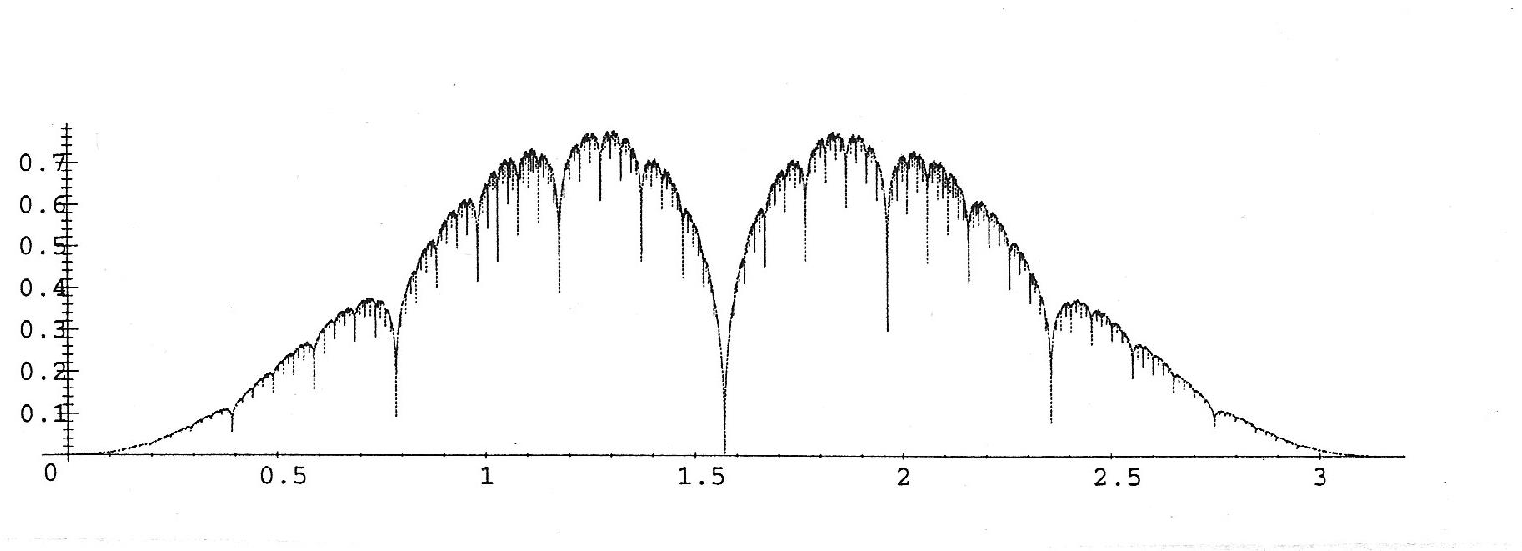}
\caption{The function}
\end{figure}

Some properties of $f(x)$ are immediately apparent.  For each factor of the infinite product, the exponent is a positive rational number with even numerator and odd denominator, so each factor is $\ge 0$ for all $x$.  Because the factors are positive powers of sine functions, they are also $\le 1$.  For each $x$, the partial products are a monotonically decreasing sequence on the interval $[0,1]$, which must approach a limiting value.  In other words, the partial products either converge to a number between 0 and 1, or they diverge to 0.  Either way, $f(x)$ is a well-defined function with values in the range $0 \le f(x) \le 1$ (in fact, $f(x)$ is strictly less than 1).

\section{Set of zeroes}

Because $\sin (2^nx)=0$ when $2^nx=m\pi$, {\it i.e.} $x=m\pi /2^n$, for any integer $m$, we have 
\begin{equation} f(m\pi/2^n)=0 \end{equation}
for every integer $m$ and non-negative integer $n$.  Thus the zeroes of $f$ are dense on the real line.

But $f(x)$ is not uniformly zero.  For example, 
\begin{equation} f(\pi/3)=(3/4)^{\pi^2/8}. \end{equation}
This follows from the fact that $2^n$ is congruent to 1, 2, or 4 mod 6, so that $\sin (2^n\pi/3)=\pm \frac 12 \sqrt 3$ and
\begin{equation} [\sin(2^n\pi/3)]^{2/(2n+1)^2}=(3/4)^{1/(2n+1)^2}.
\end{equation}  Thus
\begin{equation} f\biggl(\frac {\pi}3\biggr)=\biggl(\frac 34\biggr)^{\sum_{n=0}^{\infty}1/(2n+1)^2} \end{equation}
($\pi^2/8=\pi^2/6-\pi^2/24$, the sum of the reciprocals of all squares minus the sum for even squares.)

On the other hand, $f$ has zeroes other than $x=\pi m/2^n$.  For example, $f(x)=0$ if
\begin{equation} x=\pi \sum_{k=0}^{\infty}(-1)^k2^{-2^{2^k}}.
\end{equation}
Indeed, if $n=2^{2^j}$, then
\begin{equation} 2^nx=\pi\sum_{k=0}^{\infty}(-1)^k2^{2^{2^j}-2^{2^k}}
=\pi\sum_{k=0}^j(-1)^k2^{2^{2^j}-2^{2^k}}
+\pi\sum_{k=j+1}^{\infty}(-1)^k2^{2^{2^j}-2^{2^k}}. \end{equation}
But the sum from $k=0$ to $j$ is an integer, so
\begin{eqnarray} 0\le [\sin (2^nx)]^2=\Bigl[\sin\Bigl(
\pi\sum_{k=j+1}^{\infty}(-1)^k2^{2^{2^j}-2^{2^k}}\Bigr)\Bigr]^2\\
<\Bigl(\pi\sum_{k=j+1}^{\infty}(-1)^k2^{2^{2^j}-2^{2^k}}\Bigr)^2
=\pi^2\Bigl|\sum_{k=j+1}^{\infty}(-1)^k2^{2^{2^j}-2^{2^k}}\Bigr|^{\,2}\nonumber\\
<\pi^2\cdot 2^{(2^{2^j}-2^{2^{j+1}})2} \mbox { (since the series is alternating)}\nonumber\\ =\pi^2\cdot 2^{2(n-n^2)}. \nonumber\end{eqnarray}
It follows that
\begin{equation} 0\le[\sin(2^nx)]^{2/(2n+1)^2}<\pi^{2/(2n+1)^2}\cdot
2^{-2(n^2-n)/(2n+1)^2}. \end{equation}
Now \begin{equation} 
\lim_{n\to \infty}\frac 2{(2n+1)^2}=0,\end{equation} and
\begin{equation} \lim_{n\to \infty}\frac {-2(n^2-n)}{(2n+1)^2}=
\lim_{n\to \infty}\frac {-2n^2+2n}{4n^2+4n+1}=-\frac 12\, , \end{equation}
so \begin{equation} \lim_{n\to \infty}\pi^{2/(2n+1)^2}\cdot
2^{-2(n^2-n)/(2n+1)^2}=\pi^0\cdot 2^{-1/2}=\frac 12\sqrt 2. 
\end{equation}
Thus, for sufficiently large $n$, the upper bound in (9) gets arbitrarily close to $\frac 12\sqrt 2$, and in particular, beyond a certain point it becomes less than $\frac9{10}$, say, and stays less than $\frac 9{10}$ for all larger $n$.  A bit of experimentation reveals that this point occurs when $n=3$ (that is, $\pi^{2/49}\cdot 2^{-12/49}<\frac 9{10}$).

Thus there are an infinite number of values of $n$ (namely $n=2^{2^j}$, where $j$ is any integer $\ge 1$, so that $n>3$) for which
\begin{equation} [\sin(2^nx)]^{2/(2n+1)^2}<\frac 9{10}\, . \end{equation}
Since there are no values of $n$ for which 
\begin{equation} [\sin(2^nx)]^{2/(2n+1)^2}>1, \end{equation}
it follows that 
\begin{equation} 0\le f(x)\le \prod_{j=1}^{\infty}\frac 9{10}=
\lim_{j\to \infty}\Bigl(\frac 9{10}\Bigr)^j=0. \end{equation}
Note that $x$ is an irrational multiple of $\pi$, because the binary expansion of the sum in (6) consists of 2 zeroes, followed by 2 ones, followed by 16 zeroes, 256 ones, 65536 zeroes, etc.

Nevertheless, for ``most" $x$, $f(x)>0$.
\begin{theorem} The set of zeroes of $f(x)$ in the interval $0\le x\le \pi$ has measure $0$.
\end{theorem}

\textbf {Proof:} For each positive integer $k$, let
\begin{equation} A_k=\Bigl\{x\in [0,\pi]:\Bigl|x-\frac{m\pi}{2^n}\Bigr|
>\frac 1{2^{n+\sqrt n+k}}\mbox { for all non-negative integers $m$ and $n$}\Bigr\}. \end{equation}
Some of the intervals excluded from $A_k$ overlap, but we can obtain a lower bound on the measure of $A_k$ by subtracting from $\pi$ the lengths of all the excluded intervals.  When $m$ is even, $m/2^n$ is equal to an odd integer over a smaller power of 2, so when we add up the lengths of the excluded intervals, we can ignore even values of $m$, except for the case $m=n=0$.

Fix $n\ge 1$.  There are $2^{n-1}$ odd values of $m$ for which $m\pi/2^n$ is in the interval $[0,\pi]$, and there is an excluded interval of length $2/2^{n+\sqrt n+k}$ for each such $m$.  The total length of all these intervals is $1/2^{\sqrt n+k}$.  

Summing over all $n\ge 1$, we get 
\begin{equation} \sum_{n=1}^{\infty}\frac 1{2^{\sqrt n+k}}=\frac 1{2^k}
\sum_{n=1}^{\infty}2^{-\sqrt n}, \end{equation}
where
\begin{equation} \sum_{n=1}^{\infty}2^{-\sqrt n}<\int_0^{\infty}
2^{-\sqrt x}\,dx. \end{equation}
With the change of variables $u=(\log 2)^2x$ and $z=-\sqrt u$, we have
\begin{eqnarray} \int_0^{\infty}2^{-\sqrt x}\,dx=\int_0^{\infty}
e^{-\sqrt x\log 2}\,dx=\frac 1{(\log 2)^2}\int_0^{\infty}e^{-\sqrt u}
\,du\\ =\frac 2{(\log 2)^2}\int_0^{-\infty}ze^z\,dz=\frac 2{(\log 2)^2}
(z-1)e^z\biggr|_{\,0}^{-\infty}=\frac 2{(\log 2)^2}<\frac{25}6\,.\nonumber\end{eqnarray}
Thus
\begin{equation} \sum_{n=1}^{\infty}\frac 1{2^{\sqrt n+k}}<\frac {25}
{6\cdot 2^k}\, . \end{equation}
For $n=0$, there are two excluded intervals (around 0 and $\pi$), each with length $2^{-k}$.  So the total length of all excluded intervals is less than
\begin{equation} \frac {25}{6\cdot 2^k}+\frac 2{2^k}=\frac{37}{6\cdot
2^k}\, , \end{equation}
and the measure of $A_k$ is greater than $\pi -37/(6\cdot 2^k)$.

Next, we find a lower bound on $f(x)$ for $x\in A_k$.  Suppose $0<\delta \le \frac 12$.  Then
\begin{equation} \sin \delta >\delta -\mbox{$\frac 16$}\,\delta ^3=\delta (1-\mbox{$\frac 16$}\,\delta ^2)>\delta [1-\mbox{$\frac 16$}(\mbox{$\frac 12$})^2]=\mbox{$\frac{23}{24}$}\delta\, .\end{equation}
Now suppose $|z-\pi m|>\delta$ for every integer $m$.  Then 
\begin{equation} |\sin z|>\sin\delta >\mbox {$\frac{23}{24}$}\delta.\end{equation}
Let $n$ be any non-negative integer and suppose $|2^nx-\pi m|>\delta$.  Then $|\sin (2^nx)|>\frac{23}{24}\delta$, and the same conclusion follows from the condition
\begin{equation} \Bigl|x-\frac {\pi m}{2^n}\Bigr|>\frac{\delta}{2^n}.
\end{equation}

If $k$ is any positive integer and $n$ is any non-negative integer, then
\begin{equation} 0<\frac 1{2^{\sqrt n+k}}\le \frac 12\, ,\end{equation}
so we can let $\delta =1/2^{\sqrt n+k}$ and conclude that if
\begin{equation} \Bigl|x-\frac {\pi m}{2^n}\Bigr|>\frac 1{2^{n+\sqrt n
+k}} \end{equation}
for every integer $m$, then
\begin{equation} |\sin (2^nx)|>\frac {23}{24}\cdot \frac 1{2^{\sqrt n+k}}\, . \end{equation}
In other words, if $x\in A_k$, so that (26) holds for every integer $m$
and every non-negative integer $n$, then (27) holds for every non-negative integer $n$.  It follows that
\begin{equation} \log |\sin (2^nx)|>\log\mbox{$\frac{23}{24}$}-(\sqrt n
+k)\log 2>-\mbox{$\frac 1{23}$}-(\sqrt n+k)\log 2, \end{equation}
so \begin{eqnarray} \frac 2{(2n+1)^2}\log |\sin (2^nx)|>\frac{-2}{23(2n+1)^2}-\frac{2(\sqrt n+k)(\log 2)}{(2n+1)^2} \\ =-\frac{\frac 2{23}+2k\log 2}{(2n+1)^2}-\frac{2\sqrt n\log 2}{(2n+1)^2}.\nonumber \end{eqnarray}
Therefore
\begin{eqnarray} \log f(x)=\sum_{n=0}^{\infty}\frac 2{2n+1)^2}\log
|\sin (2^nx)| \\ >-\sum_{n=0}^{\infty}\frac{\frac 2{23}+2k\log 2}{(2n+1)^2}-(\log 2)\sum_{n=0}^{\infty}\frac{2\sqrt n}{(2n+1)^2}\, ,\nonumber
\end{eqnarray} where
\begin{equation} \sum_{n=0}^{\infty}\frac{\frac 2{23}+2k\log 2}{(2n+1)^2}=\frac {\pi^2}8\Bigl(\frac 2{23}+2k\log 2\Bigr) \end{equation} and
\begin{eqnarray} \sum_{n=0}^{\infty}\frac{2\sqrt n}{(2n+1)^2}=
\sum_{n=1}^{\infty}\frac{2\sqrt n}{(2n+1)^2}<\sum_{n=1}^{\infty}\frac{2\sqrt n}{(2n)^2}=\sum_{n=1}^{\infty}\frac 1{2n^{3/2}} \\
=\frac 12+\sum_{n=2}^{\infty}\frac 1{2n^{3/2}}<\frac 12+
\int_1^{\infty}\frac 1{2x^{3/2}}\,dx=\frac 12+1=\frac 32\, .\nonumber
\end{eqnarray} Thus
\begin{equation} \log f(x)>\frac{-\pi^2}8\Bigl(\frac 2{23}+2k\log 2\Bigr)-\frac 32\log 2>-1.147-1.7103k \end{equation} and
\begin{equation} f(x)>e^{-1.147-1.17103k}>\frac 1{(3.15)(5.531)^k}\, .
\end{equation}
It follows that if $f(x)\le 1/(3.15)(5.531)^k$, then 
$x\not {\! \! \in}A_k$.
Therefore, for every positive integer $k$, the measure of the set of all $x$ in $[0,\pi]$ for which $f(x)\le 1/(3.15)(5.531)^k$ is less than $37/(6\cdot 2^k)$, the measure of the complement of $A_k$.  The set of $x$ for which $f(x)=0$ is a subset of the set of $x$ for which 
$f(x)\le 1/(3.15)(5.531)^k$ for every positive integer $k$.  Therefore the measure of the set of $x$ for which $f(x)=0$ is less than $37/(6\cdot 2^k)$ for every $k$, and is thus 0.

\section{Points of continuity}

A function $f(x)$ is said to be upper semicontinuous \cite{gelbaum} (p. 22) at $x=a$ if for every $\varepsilon >0$, there exists $\delta >0$ such that $f(x)<f(a)+\varepsilon$ whenever $|x-a|<\delta$.  Note the asymmetry of this definition: $f(x)$ must be less than $f(a)+\varepsilon$ but need not be greater than $f(a)-\varepsilon$.  Note also that continuity implies upper semicontinuity.

We shall prove that our function $f(x)$ is upper semicontinuous at all $x$.  Two corollaries are that $f(x)$ is continuous at $x$ if and only if $f(x)=0$, and that $f(x)$ is Lebesgue integrable.

\begin{theorem} $f(x)$ is upper semicontinuous at all $x$. \end{theorem}

\textbf{Proof:} Fix $x$.  We consider two cases.

Case 1: $\sin (2^nx)=0$ for some $n$.  In this case, $2^nx=m\pi$ for some integer $m$, so $x=\pi m/2^n$.  Then, because
\begin{equation} 0\le [\sin (2^jt)]^{2/(2j+1)^2}\le 1 \end{equation}
for all $t$ and all $j$, we have
\begin{equation} 0\le f(t)\le [\sin(2^nt)]^{2/(2n+1)^2} \end{equation}
for all $t$.  Since the right side of (36) is a continuous function of $t$ which is equal to 0 at $t=x$, we have $\lim_{t\to x}f(t)=f(x)$ by the squeeze theorem.  Therefore $f(t)$ is continuous at $t=x$, and thus upper semicontinuous.

Case 2: $\sin (2^nx)$ is not equal to 0 for any non-negative integer $n$.  In this case, neither is $\sin (2^nx)$ equal to $\pm 1$ for any $n$.  For if $\sin (2^nx)=\pm 1$, then $2^nx$ is an odd multiple of $\pi /2$, in which case $2^{n+1}x$ is a multiple of $\pi$, and $\sin (2^{n+1}x)=0$.  Thus
\begin{equation} 0<[\sin(2^nx)]^{2/(2n+1)^2}<1 \end{equation}
for all $n$.  It follows that the partial products
\begin{equation} f_k(x)=\prod_{n=0}^k[\sin(2^nx)]^{2/(2n+1)^2} \end{equation}
decrease monotonically with $k$, approaching $f(x)$ as $k\to\infty$.

Given $\varepsilon >0$, we must find $\delta >0$ such that $f(t)<f(x)+\varepsilon $ whenever $|t-x|<\delta $.  Let $k$ be the least non-negative integer such that $f_k(x)<f(x)+\varepsilon ^2$.  Note that $k$
must exist, because $f_k(x)$ approaches $f(x)$ monotonically from above as $k\to \infty$.  Let $\lambda =f_k(x)$.  Note that $k$ and 
$\lambda$ depend only on $x$ and $\varepsilon$.  

If $\varepsilon \le \frac 17$, let 
\begin{equation} \delta =\frac {(2k+1)^2\lambda ^{2k+1)^2/2}\varepsilon}{21\cdot 2^{k+1}} \, .\end{equation}
If $\varepsilon>\frac 17$, let $\delta$ have the same value that it would have for $\varepsilon =\frac 17$.  Note that because $\delta$ is a function of $k$, $\lambda$, and $\varepsilon$, $\delta$ also depends only on $x$ and $\varepsilon$.  We will show that if $|t-x|<\delta$, then $f(t)<f(x)+\varepsilon$.

\begin{lemma} If $\varepsilon \le \frac 17$, then $\lambda <\frac {81}{100}$.\end{lemma}

\textbf{Proof:} Suppose $k=0$.  Then
\begin{eqnarray} f_k(x)-f(x)=f_0(x)-f(x)>f_0(x)-f_1(x)\\
=(\sin^2x)[1-(\sin 2x)^{2/9}]=(\sin^2x)[1-(2\sin x\cos x)^{2/9}]\nonumber \\ \ge (\sin^2x)[1-(2\cos x)^{2/9}]=(\sin^2x)[1-(4\cos^2x)^{1/9}]\nonumber \\ (\sin^2x)[1-4^{1/9}(1-\sin^2x)^{1/9}]. \nonumber
\end{eqnarray}
Now suppose $\sin x\ge \frac 9{10}$.  Then $\sin^2x\ge \frac{81}{100}$,
$1-\sin^2x\le \frac{19}{100}$, 
\begin{equation} 4(1-\sin^2x)\le \mbox{$\frac{76}{100}$,}\end{equation}
\begin{equation} 4^{1/9}(1-\sin^2x)^{1/9}\le (\mbox{$\frac{76}{100}$})^{1/9}<\mbox{$\frac{97}{100}$,}\end{equation}
\begin{equation} 1-4^{1/9}(1-\sin^2x)^{1/9}>\mbox{$\frac 3{100}$,}
\end{equation} and
\begin{equation} (\sin^2x)[1-4^{1/9}(1-\sin^2x)^{1/9}]>(\mbox{$\frac
{81}{100}$})(\mbox{$\frac 3{100}$})>\mbox{$\frac 1{49}$}\ge \varepsilon^2.\end{equation}
But by our definition of $k$, $f_k(x)-f(x)<\varepsilon^2$.  This contradiction establishes that if $k=0$, then $\sin x<\frac 9{10}$.
Therefore
\begin{equation} \lambda =f_0(x)=\sin^2x<(\mbox{$\frac 9{10}$})^2=
\mbox{$\frac{81}{100}$.} \end{equation}

Suppose $k=1$.  Then 
\begin{equation} f_k(x)=f_1(x)=(\sin x)^2(\sin 2x)^{2/9},\end{equation}
which has its maximum value, just below $\frac{81}{100}$, when $\sin x
=\sqrt {10/11}$.  If $k>1$, then the maximum value of $f_k(x)$ is less than the maximum value of $f_1(x)$.

\begin{lemma} For every non-negative integer $n$, and every real number $x$ and $t$ such that $t\ne x$, and such that $\sin(2^nx)$ and 
$\sin(2^nt)$ are not $0$, we have
\begin{equation} \log [\sin(2^nt)]^{2/(2n+1)^2}-\log [\sin(2^nx)]^{2/(2n+1)^2}<\frac{2^{n+1}|t-x|}{(2n+1)^2|\sin(2^nx)|}\, .
\end{equation}
\end{lemma}

\textbf{Proof:} Since the righthand side of the above inequality is always positive, the inequality is satisfied whenever the lefthand side is negative.  Suppose the lefthand side is positive.

Let \begin{equation} g_n(t)=\log[\sin(2^nt)]^{2/(2n+1)^2}.
\end{equation}
Then $g_n(t)$ is periodic with period $2^{-n}\pi$, and because $g_n(t)$ is an even function, we have
\begin{equation} g_n(t)=g_n(t+2^{-n}j\pi)=g_n(-t+2^{-n}j\pi)
\end{equation}
for all integers $j$.  Furthermore, $g_n(t)$ increases monotonically on each interval
\begin{equation} \frac{m\pi}{2^{n+1}}<t<\frac{(m+1)\pi}{2^{n+1}},
\end{equation}
increasing if $m$ is even, and decreasing if $m$ is odd.  Indeed, if 
$g_n(t^{\prime})=g_n(t)$, then $t^{\prime}$ must be equal to $t+2^{-n}j\pi$ or $-t+2^{-n}j\pi$ for some integer $j$.

Suppose $x$ and $t$ are both between $m\pi/2^{n+1}$ and $(m+1)\pi/2^{n+1}$ for some integer $m$, and suppose that $t^{\prime}$ is a real number, not between $m\pi/2^{n+1}$ and $(m+1)\pi/2^{n+1}$, such that 
$g_n(t^{\prime})=g_n(t)$.  Then $t^{\prime}=t+2^{-n}j\pi$ or $-t+2^{-n}j\pi$ for some integer $j$.  If $t^{\prime}=t+2^{-n}j\pi$, then $j\ne 0$, so $|t^{\prime}-t|\ge 2^{-n}\pi$, but $|t-x|<\pi/2^{n+1}$, so 
$|t^{\prime}-x|>\pi/2^{n+1}$ and $|t^{\prime}-x|>|t-x|$. If If $t^{\prime}=-t+2^{-n}j\pi$, let If $x^{\prime}=-x+2^{-n}j\pi$.  Then $|t^{\prime}-x^{\prime}|=|t-x|$, and either $x<t<t^{\prime}<x^{\prime}$, 
$x^{\prime}<t^{\prime}<t<x$, $t<x<x^{\prime}<t^{\prime}$, or $t^{\prime}<x^{\prime}<x<t$.  Whichever case applies, we have $|t^{\prime}-x|>|t-x|$.  So if
\begin{equation} g_n(t)-g_n(x)<\frac{2^{n+1}|t-x|}{(2n+1)^2|\sin(2^nx)|},\end{equation}
then $g_n(t^{\prime})-g_n(x)$, which is equal to $g_n(t)-g_n(x)$, is less than
\begin{equation}\frac{2^{n+1}|t^{\prime}-x|}{(2n+1)^2|\sin(2^nx)|}.\end{equation}
Thus it suffices to prove the lemma for the case where $x$ and $t$ are both between $m\pi/2^{n+1}$ and $(m+1)\pi/2^{n+1}$ for the same integer $m$.

Now \begin{equation} [\sin(2^nt)]^{2/(2n+1)^2}\ge 0 \end{equation}
for all $t$, so
\begin{equation} g_n(t)=\log|\sin(2^nt)|^{2/(2n+1)^2}=\frac 2{(2n+1)^2}
\log|\sin(2^nt)|. \end{equation}
We want to find an upper bound on $|g_n^{\prime}(t)|$.  We have 
\begin{equation} |\sin u|=\left \{ \begin{array}{ll}
\sin u & \mbox{if }\sin u>0\\ -\sin u & \mbox{if }\sin u<0 \end{array}
\right \} \end{equation} so
\begin{equation} \frac d{du}|\sin u|=\left \{ \begin{array}{ll}
\cos u & \mbox{if }\sin u>0\\ -\cos u & \mbox{if }\sin u<0 \end{array}
\right \} \end{equation} and
\begin{equation} \frac d{du}\log |\sin u|=\cot u\,\mbox{ if }\sin u\ne 0.\end{equation} Thus
\begin{equation} g_n^{\prime}(t)=\frac{2^{n+1}}{(2n+1)^2}\cot(2^nt)
\end{equation} and
\begin{equation} |g_n^{\prime}(t)|=\frac{2^{n+1}}{(2n+1)^2}|\cot(2^nt)|
<\frac{2^{n+1}}{(2n+1)^2}\cdot \frac 1{|\sin (2^nt)|}\, .\end{equation}
By hypothesis, the expression on the lefthand side of (47) is positive; that is
\begin{equation} \log [\sin(2^nt)]^{2/(2n+1)^2}>\log [\sin(2^nx)]^{2/(2n+1)^2} \end{equation}
so $|\sin(2^nt)|>|\sin(2^nx)|$  and
\begin{equation} |g_n^{\prime}(t)|<\frac{2^{n+1}}{(2n+1)^2|\sin(2^nx)|}
\, .\end{equation}
The same upper bound holds for $|g_n^{\prime}(u)|$ when $u$ is between $x$ and $t$, and we also know that $g_n^{\prime}(u)$ has the same sign over this entire interval, because $x$ and $t$ are both between $m\pi/2^{n+1}$ and $(m+1)\pi/2^{n+1}$ for the same integer $m$.
It follows that
\begin{eqnarray} g_n(t)-g_n(x)=|g_n(t)-g_n(x)|=\left|\int_x^t
g_n^{\prime}(u)\, du\right| \\ =\left|\int_x^t|g_n^{\prime}(u)|\, du
\right|<\frac{2^{n+1}|t-x|}{(2n+1)^2|\sin(2^nx)|}\, .\nonumber \end{eqnarray}
This concludes the proof of Lemma 2.

$\,\,$

We now continue the proof of Theorem 2. Defining $g_n$ as in (48), and assuming $x$ is not a multiple of $2^{-k}\pi$, so that $f_k(x)\ne 0$, we have
\begin{equation} \log f_k(x)=\sum_{n=0}^kg_n(x), \end{equation}
and the same holds for $t$.  It follows from Lemma 2 that
\begin{equation} \log f_k(t)-\log f_k(x)<\sum_{n=0}^k\frac{2^{n+1}|t-x|}{(2n+1)^2|\sin(2^nx)|}=|t-x|\sum_{n=0}^k\frac{2^{n+1}}{(2n+1)^2|\sin(2^nx)|}\, .\end{equation}
From (37) we have
\begin{equation} \lambda =f_k(x)=\prod_{n=0}^k[\sin(2^nx)]^{2/(2n+1)^2}
<[\sin(2^nx)]^{2/(2n+1)^2}=|\sin(2^nx)|^{2/(2n+1)^2}\end{equation}
for each $n$, so
\begin{equation} |\sin(2^nx)|>\lambda^{(2n+1)^2/2} \end{equation} and
\begin{equation} \log f_k(t)-\log f_k(x)<|t-x|\sum_{n=0}^k\frac{2^{n+1}}{(2n+1)^2\lambda^{(2n+1)^2/2}}\, .\end{equation}

Let \begin{equation} a_n=\frac{2^{n+1}}{(2n+1)^2\lambda^{(2n+1)^2/2}}\, .\end{equation} Then
\begin{equation}\frac{a_{n+1}}{a_n}=\frac{2(2n+1)^2}{(2n+3)^2\lambda^{4n+4}}\, .\end{equation}
By Lemma 1, $\lambda<\frac{81}{100}$, so if $n\ge 1$, then
\begin{equation} \frac{a_{n+1}}{a_n}>\frac 29\left(\frac{81}{100}\right)^{-8}>\frac {18}5 \end{equation}
and if $n=0$ then
\begin{equation} \frac{a_{n+1}}{a_n}>\frac 29\left(\frac{81}{100}\right)^{-4}>\frac 12\, . \end{equation}
It follows that
\begin{equation} \sum_{n=0}^k\frac{2^{n+1}}{(2n+1)^2\lambda^{(2n+1)^2/2}}<\frac{3\cdot 2^{k+1}}{(2k+1)^2\lambda^{(2k+1)^2/2}}\, ,
\end{equation} and
\begin{equation} \log f_k(t)-\log f_k(x)<\frac{3\cdot 2^{k+1}|t-x|}
{(2k+1)^2\lambda^{(2k+1)^2/2}}\, .\end{equation}

We already know that $f(t)<f_k(t)$, because $f_j(t)$ decreases monotonically with $j$, approaching $f(t)$ in the limit as $j\to \infty$.  We also know that $f(t)>0$, because $f(t)>f(x)$ by hypothesis, so $\log f(t)$ is defined.  Thus $\log f(t)<\log f_k(t)$ and \begin{equation} \log f(t)<\log f_k(x)+\frac{3\cdot 2^{k+1}|t-x|}
{(2k+1)^2\lambda^{(2k+1)^2/2}}\, .\end{equation}
We also have, by the definition of $k$, that $f_k(x)<f(x)+\varepsilon^2$.  We must show that for $|t-x|$ sufficiently small, we have
\begin{equation} \log[f(x)+\varepsilon^2]+\frac{3\cdot 2^{k+1}|t-x|}
{(2k+1)^2\lambda^{(2k+1)^2/2}}<\log[f(x)+\varepsilon]. \end{equation}

To prove this, we consider two cases.

Case 1: $f(x)<\frac 14\varepsilon$.  Then, if $f(x)\le\frac 17$, we have $\varepsilon^2\le\frac 17\varepsilon$, so $f(x)+\varepsilon^2<
\frac 14\varepsilon +\frac 17\varepsilon <\frac 12\varepsilon$, and
\begin{equation} \log[f(x)+\varepsilon^2]<\log(\mbox{$\frac 12$}\varepsilon)=\log \varepsilon -\log 2. \end{equation}  But
\begin{equation} \log[f(x)+\varepsilon]\ge \log \varepsilon 
\end{equation} so it suffices to have
\begin{equation} \frac{3\cdot 2^{k+1}\delta}{(2k+1)^2\lambda^{(2k+1)^2/2}}<\log 2, \end{equation} or
\begin{equation} \delta <\frac{(\log 2)(2k+1)^2\lambda^{(2k+1)^2/2}}
{3\cdot 2^{k+1}}\, .\end{equation}

Case 2: $f(x)\ge \frac 14\varepsilon$.  Then
\begin{eqnarray} \log [f(x)+\varepsilon^2]=\log \left\{f(x)\left[1+
\frac{\varepsilon^2}{f(x)}\right]\right\} \\ =\log f(x)+\log \left[
1+\frac{\varepsilon^2}{f(x)}\right ]<\log f(x)+\frac{\varepsilon^2}{f(x)}\, .\nonumber\end{eqnarray} Also
\begin{equation} \log [f(x)+\varepsilon]=\log f(x)+\log \left[
1+\frac{\varepsilon}{f(x)}\right ]. \end{equation}
Since $\varepsilon /f(x)\le 4$, and the $\log$ function is concave down, we have
\begin{equation} \log \left[1+\frac{\varepsilon}{f(x)}\right ]\ge
\frac{\log 5}4\cdot\frac{\varepsilon}{f(x)}\, .\end{equation}
Therefore
\begin{equation}\ \ \ \ \ \log [f(x)+\varepsilon]-\log [f(x)+\varepsilon^2]\ge \frac{\log 5}4\cdot\frac{\varepsilon}{f(x)}-\frac{\varepsilon^2}{f(x)}=\left[\frac{\log 5}4-\varepsilon\right]\frac{\varepsilon}{f(x)}\, ,\end{equation} which is 
\begin{equation} \ge \left[\frac{\log 5}4-\frac 17\right]\,\frac{\varepsilon}{f(x)}>\frac 17\frac{\varepsilon}{f(x)}>\frac 17\,\varepsilon\,\mbox{ if }\,\varepsilon\le\frac 17\, .\end{equation}
So if
\begin{equation} \frac{3\cdot 2^{k+1}\delta}{(2k+1)^2\lambda^{(2k+1)^2/2}}<\frac 17\varepsilon, \end{equation}
{\it i.e.} if
\begin{equation} \delta<\frac{(2k+1)^2\lambda^{(2k+1)^2/2}\varepsilon}
{21\cdot 2^{k+1}}\, ,\end{equation} then
\begin{equation} \log f(t)<\log[f(x)+\varepsilon] \end{equation}
and \begin{equation} f(t)<f(x)+\varepsilon . \end{equation}

If $\varepsilon \le \frac 17$, then $\frac 17\varepsilon\le\frac 1{49}
<\log 2$, so whether or not $f(x)<\frac 14\varepsilon$, if
\begin{equation} |t-x|<\frac{(2k+1)^2\lambda^{(2k+1)^2/2}\varepsilon}
{21\cdot 2^{k+1}}\, ,\end{equation}
then $f(t)<f(x)+\varepsilon$.  This concludes the proof of Theorem 2.

$\,\,$

\textbf{Corollary:} {\it $f(x)$ is continuous at $x$ if and only if
$f(x)=0$.}

$\,\,$

\textbf{Proof:} Because the set of zeroes of $f(x)$ is everywhere dense, $f(x)$ cannot be continuous if $f(x)\ne 0$.  On the other hand, $f(x)$ is upper semicontinuous everywhere, so given $x$, for every $\varepsilon >0$, there exists $\delta >0$ such that if $|t-x|<\delta$, then $f(t)<f(x)+\varepsilon$.  But $f(t)$ is never negative, so if $f(x)=0$, then $f(t)>f(x)-\varepsilon$.  Therefore $f(x)$ is continuous at $x$ if $f(x)=0$.

$\,\,$

Another corollary of Theorem 2 is that $f(x)$ is Lebesgue integrable, because a function that is bounded from below and upper semicontinuous on a closed interval is Lebesgue integrable over that interval \cite {hawkins} (p. 151).  Indeed, suppose that $f(x)$ is upper semicontinuous on $[a,b]$, and let $r$ be a lower bound.  Let 
$s$ be an upper bound of $f(x)$, which must exist, because if $\{x_i\}$ is a sequence of real numbers on which $f$ is unbounded, then $f$ cannot be upper semicontinuous on an accumulation point of $\{x_i\}$.  Now suppose $f(x)<y$ for some $x\in [a,b]$ and $y\in [r,s]$.  Let $\varepsilon =y-f(x)$.  There exists $\delta >0$ such that $f(t)<f(x)+\varepsilon =y$ whenever $|t-x|<\delta$.  In other words, if $f(x)<y$, then there is a neighborhood $U$ of $x$ such that $f(t)<y$ for all $t$ in $U$.  It follows that for every $y$ in $[r,s]$, the set
\begin{equation} S_y=\{x\in [a,b]:f(x)<y\} \end{equation}
is an open set of $[a,b]$.  Let $g(y)$ be the measure of $S_y$.  Then $g(y)$ increases monotonically on the interval $[r,s]$, so $g(y)$ is Riemann integrable.  But the Lebesgue integral $\int_a^bf(x)\,dx$ is equal to the Riemann integral $\int_r^sg(y)\,dy$.

\section{A lower bound on the Lebesgue integral}

\begin{theorem} The Lebesgue integral $\int_0^{\pi}f(x)\,dx$ is strictly positive. \end{theorem}

\textbf{Proof:} First, we prove that for all $k$, the improper integral
\begin{equation} \int_0^{\pi}\log f_k(x)\,dx \end{equation}
converges to a value $>-\pi^3/4$.  We have
\begin{equation} \log f_k(x)=\sum_{n=0}^k\log [\sin(2^nx)]^{2/(2n+1)^2}
=\sum_{n=0}^k\frac 2{(2n+1)^2}\log |\sin(2^nx)| \end{equation}
(where we take both sides to be $-\infty$ when $x$ is a multiple of
$2^{-k}\pi$), so
\begin{equation} \int_0^{\pi}\log f_k(x)\,dx=\sum_{n=0}^k\frac 2{(2n+1)^2}\int_0^{\pi}\log |\sin(2^nx)|\,dx,\end{equation}
where the integrals on both sides are improper.  By the change of variables $u=2^nx$, we have 
\begin{equation} \int_0^{\pi}\log |\sin(2^nx)|\,dx=2^{-n}\int_0^{2^n\pi}\log |\sin u|\,du, \end{equation}
which, by the symmetry and periodicity of the sine function, is equal to \begin{equation} 2\int_0^{\pi/2}\log |\sin u|\,du. \end{equation}
For $0\le u\le \pi/2$, we have $|\sin u|=\sin u\ge 2u/\pi$, so
\begin{equation} \log |\sin u|\ge \log\left(\frac{2u}{\pi}\right),
\end{equation}
and, by the change of variables $z=2u/\pi$, we have
\begin{equation} 2\int_0^{\pi/2}\log |\sin u|\,du\ge 2\int_0^{\pi/2}
\log\left(\frac{2u}{\pi}\right)\,du=\pi\int_0^1\log z\,dz=-\pi.
\end{equation}
Therefore
\begin{equation} \int_0^{\pi}\log f_k(x)\,dx\ge\sum_{n=0}^k\frac{-2\pi}
{(2n+1)^2}>\sum_{n=0}^\infty\frac{-2\pi}{(2n+1)^2}=-2\pi\cdot\frac{\pi^2}8=-\frac{\pi^3}4\, .\end{equation}

Now, for each positive integer $k$, let
\begin{equation} B_k=\{ x\in[0,\pi]:f_k(x)>e^{-\pi^2/2}\}. 
\end{equation}
Because $f_k(x)$ is continuous, $B_k$ is open, and hence measurable.  We want to show that the measure of $B_k$ is $>\pi/2$ for all $k$.  The complement of $B_k$ in $[0,\pi]$ is
\begin{equation}\ \ \ \ \ \bar B_k=\{x\in[0,\pi]:f_k(x)\le e^{-\pi^2/2}\}=
\{x\in[0,\pi]:\log f_k(x)\le -\pi^2/2\}, \end{equation}
a closed set.  Suppose the measure of $\bar B_k$ is $\ge \pi/2$ for some $k$.  Then, because $\log f_k(x)\le 0$ (incl. $-\infty$) for all $x$, and $\bar B_k$ is a subset of $[0,\pi]$, we have
\begin{equation} \int_0^{\pi}\log f_k(x)\,dx\le\int_{\bar B_k}\log f_k(x)\,dx\le -\frac{\pi^2}2\cdot\frac{\pi}2=-\frac{\pi^3}4\, ,
\end{equation}
but by (98), this integral is $>-\pi^3/4$, and this contradiction establishes that the measure of $\bar B_k$ is $<\pi/2$, and the measure of $B_k$ is $>\pi/2$.

Since $f_k(x)\ge 0$ for all $x$, we have
\begin{equation} \int_0^{\pi}f_k(x)\,dx\ge \int_{B_k}f_k(x)\,dx>\frac{\pi}2\,e^{-\pi^2/2} \end{equation}
for all $k$.

Now the sequence $\{f_k(x)\}_{k=0}^{\infty}$ converges pointwise to $f(x)$ on the interval $0\le x\le \pi$, and $0\le f_k(x)\le 1$ for all $x$ and all $k$.  It follows from the Lebesgue dominated convergence theorem \cite {bressoud} (p. 183, Theorem 6.19) that
\begin{equation} \int_0^{\pi}f(x)\,dx=\lim_{k\to \infty}\int_0^{\pi}
f_k(x)\,dx. \end{equation}
But $\int_0^{\pi}f_k(x)\,dx$ exists and is greater than $\frac{\pi}2
e^{-\pi^2/2}$ for each $k$.  It follows that $\int_0^{\pi}f(x)\,dx$ exists and is $\ge\frac{\pi}2e^{-\pi^2/2}$.

This concludes the proof of Theorem 3.

$\,\,$

An immediate consequence of Theorem 3 is that $f(x)$ is not Riemann integrable.  If it were, then the Riemann integral would be equal to the Lebesgue integral, but because the zeroes of $f(x)$ are dense on the interval $[0,\pi]$, every lower Riemann sum is zero.

However, we do have

\begin{theorem} The lim inf of the upper Riemann sums of $f(x)$ is equal to the Lebesgue integral. \end{theorem}

\textbf{Proof:} Let $R$ be the set of partitions of the interval $[0,\pi]$ into a finite number of intervals, and let $L$ be the set of partitions of $[0,\pi]$ into a finite number of Borel sets.  Let
\begin{equation} U_R=\liminf_{{\cal P}\in R}\sum_{S\in{\cal P}}\mu(S)
\limsup_{x\in S}f(x) \end{equation} and let
\begin{equation} U_L=\liminf_{{\cal P}\in L}\sum_{S\in{\cal P}}\mu(S)
\limsup_{x\in S}f(x) \end{equation}
where $\mu$ is Borel measure.  Because every interval is a Borel set, $R$ is a subset of $L$, and $U_R\ge U_L$.

For each $k$, let
\begin{equation} U_{R,k}=\liminf_{{\cal P}\in R}\sum_{S\in{\cal P}}\mu(S)\limsup_{x\in S}f_k(x). \end{equation}
For each $k$, $f(x)\le f_k(x)$ for all $x$, so $U_R\le U_{R,k}$.  Also, for each $k$, $f_k(x)$ is continuous, and hence Riemann integrable.  By the Lebesgue dominated convergence theorem \cite {bressoud} (p. 183), the Riemann integral of $f(x)$, and hence $U_L$, is equal to the limit as $k\to \infty$ of the Lebesgue integral of $f_k(x)$ (and hence the Riemann integral of $f_k(x)$, and $U_{R,k}$).  Therefore $U_R\le U_L$.  Since $U_R$ is both $\ge$ and $\le U_L$, $U_R=U_L$, and since $f(x)$ is Lebesgue integrable, $U_R$ is equal to the Lebesgue integral.

\section{A numerical estimate}

How can we find a decimal value for $\int _0^\pi f(x)\,dx$?  The usual numerical integration methods, such as Simpson's rule, are unstable for this function.  However, $\int _0^\pi f_k(x)\,dx$ converges to $\int _0^\pi f(x)\,dx$ as $k\to \infty$, and $f_k(x)$ is continuous, so we can estimate $\int _0^\pi f(x)\,dx$ by estimating $\int _0^\pi f_k(x)\,dx$.

Let $M_k$ be the midpoint estimate of $2\int_0^{\pi/2}f_{k+1}(x)\,dx$
(which is equal to $\int _0^\pi f_{k+1}(x)\,dx$) with $2^k$ intervals.
Then
\begin{equation} M_k=\frac {\pi}{2^k}\sum_{j=1}^{2^k}\bigl [\sin \bigl ( (2j-1)2^{n-k-2}\pi \bigr ) \bigr ]^{2/(2n+1)^2}. \end{equation}
(In the above equation, we need only compute the product up to $n=k$, instead of $n=k+1$, because the sine of any odd multiple of $\frac {\pi}2$ is $1$.)  Let
\begin{equation} M_{\infty}=\lim_{k\to\infty}M_k.\end{equation}
We conjecture that $M_{\infty}$ exists and is equal to $\int _0^\pi f(x)\,dx$.  This does not, of course, follow from the fact that for fixed $k$, the midpoint estimate of $\int _0^\pi f_{k+1}(x)\,dx$ with $2^m$ intervals converges to this integral as $m\to\infty$, and $\int _0^\pi f_{k+1}(x)\,dx$ converges to $\int _0^\pi f(x)\,dx$ as $k\to\infty$.  

Table 1 shows the values of $M_k$ for $6\le k\le 29$, in column 2.  Column 3 shows the reciprocal square roots of the differences $M_{k-1}-M_k$.  The fact that these grow linearly with $k$ means that the differences decrease as $1/k^2$, which is what we would expect, given that $\int _0^{\pi}[\sin (2^nx)]^{2/(2n+2)^2}\,dx=\pi -O(1/n^2)$.  This in turn suggests that the errors $M_k-M_{\infty}$ decrease as $1/k$.  We might expect that for a suitable choice of constants $a$ and $b$, $M_k-a(k-b)^{-1}$ should converge to $M_{\infty}$ much more rapidly than $M_k$ itself.  A bit of trial and error reveals that the values $a=.4044$ and $b=.27$ work nicely.  Column 4 shows the values of $M_k-.4044(k-.27)^{-1}$.  To 5 decimal places, $\int _0^\pi f(x)\,dx$ appears to be $1.16993\dots $

\begin{table}[htb]
\begin{tabular}{r|c|r|c}
$k$ & $M_k$ &$(M_{k-1}-M_k)^{-1/2}$& $M_k-.4044/(k-.27)$\\ \hline
6 & 1.2419727451 &   & 1.1713968 \\
7 & 1.2311527243 & 9.613598 & 1.1710636 \\
8 & 1.2230892609 & 11.136255 & 1.1707736 \\
9 & 1.2168748353 & 12.685264 & 1.1705518 \\
10 & 1.2119511226 & 14.251272 & 1.1703889 \\
11 & 1.2079596568 & 15.828283 & 1.1702709 \\
12 & 1.2046613111 & 17.412130 & 1.1701856 \\
13 & 1.2018911808 & 18.999838 & 1.1701237 \\
14 & 1.1995322446 & 20.589315 & 1.1700785 \\
15 & 1.1974993737 & 22.179160 & 1.1700452 \\
16 & 1.1957292786 & 23.768496 & 1.1700204 \\
17 & 1.1941739924 & 25.356823 & 1.1700019 \\
18 & 1.1927965318 & 26.943897 & 1.1699877 \\
19 & 1.1915679404 & 28.529638 & 1.1699769 \\
20 & 1.1904652307 & 30.114067 & 1.1699685 \\
21 & 1.1894699246 & 31.697256 & 1.1699620 \\
22 & 1.1885669999 & 33.279304 & 1.1699568 \\
23 & 1.1877441184 & 34.860317 & 1.1699527 \\
24 & 1.1869910513 & 36.440403 & 1.1699493 \\
25 & 1.1862992466 & 38.019661 & 1.1699466 \\
26 & 1.1856614980 & 39.598181 & 1.1699444 \\
27 & 1.1850716898 & 41.176044 & 1.1699426 \\
28 & 1.1845245979 & 42.753322 & 1.1699411 \\
29 & 1.1840157324 & 44.330078 & 1.1699399 
\end{tabular}
\caption{}
\end{table}

I would like to thank Daniel Asimov for his assistance computing the numbers in Table 1.


\begin{thebibliography}{3}

\bibitem{bressoud} David M. Bressoud, A Radical Approach to Lebesgue's Theory of Integration, Cambridge U. Press, 2008.

\bibitem{gelbaum} Bernard R. Gelbaum and John M. H. Olmsted, Counterexamples in Analysis, Holden Day, 1964.

\bibitem{hawkins} Thomas Hawkins, Lebesgue's Theory of Integration: Its Origins and Development, Am. Math. Soc., 2001.

\end{thebibliography}
\end{document}